\documentclass{amsart}
\usepackage[top=2cm,left=2cm,bottom=2cm,right=2cm]{geometry}
\usepackage{amssymb,amsmath,amsthm,latexsym,verbatim}

\newcounter{thmcount}
\newtheorem{theorem}[thmcount]{Theorem}

\newtheorem{definition}[thmcount]{Definition}
\newtheorem{lemma}[thmcount]{Lemma}
\newtheorem*{remark}{Remark}
\newtheorem*{theorem*}{Theorem}

\begin{document}

\bibliographystyle{amsalpha}

\title{Superstability and Finite Time Extinction for $C_{0}$-Semigroups}

\author{D. Creutz}
\address{Department of Mathematics, Vanderbilt University, Nashville, USA}
\email{darren.creutz@vanderbilt.edu}

\author{M. Mazo, Jr.}
\address{Delft Center for Systems and Control, Delft University of Technology, The Netherlands}
\email{m.mazo@tudelft.nl}

\author{C. Preda}
\address{Department of Economics and Business Modeling, University of West Timisoara, Romania}
\email{ciprian.preda@feaa.uvt.ro}

\date{12 June 2008}

\begin{abstract}
A new approach to superstability and finite time extinction of strongly continuous semigroups is presented, unifying known results and providing new criteria for these conditions to hold analogous to the well-known Pazy condition for stability.  That finite time extinction implies superstability which is in turn equivalent to several (both known and new) conditions follow from this new approach in a consistent fashion.  Examples showing that the converse statements fail are constructed, in particular, an answer to a question of Balakrishnan on superstable systems not exhibiting finite time extinction.
\end{abstract}

\maketitle

{\allowdisplaybreaks
\large
\section{Introduction}
The study of $C_{0}$-semigroups as a means to understand systems, particularly systems modelled by (partial) differential equations, has a long and rich history.  
In the control of partial differential equations (PDE) the theory of $C_{0}$-semigroups plays a fundamental role allowing to apply similar ideas as those employed in the case of ordinary (finite) differential equations.

Consider the strongly continuous semigroup $T(t)$ with generator $A$, arising from e.g. a PDE model. One can associate to this operator an alternative description in the form of a linear continuous-time system $\Sigma(A,B,C)$:
\begin{eqnarray*}
x'(t)&=&Ax(t)+Bu(t)\\
x(0)&=&0\\
y(t)&=&Cx(t),
\end{eqnarray*}
where $B$ and $C$ are bounded linear operators defined on a suitable function space. One can also define the transfer function of such a system, i.e.:  an analytic function $G$ bounded on the open right-hand complex half-plane that satisfies the condition
$$\hat{y}(\lambda)=G(\lambda)\hat{u}(\lambda),\;\lambda\in \mathbb{C}_+.$$
It is well-known that the transfer function is given by 
$$G(\lambda)= C(\lambda I-A)^{-1}B,\; \lambda\in\mathbb{C},\; Re(\lambda)>\omega_0$$
where $\omega_0$ is the growth bound of the semigroup generated by $A$.
The representation of these systems by transfer functions is the key to more sophisticated control theory.

A central concept in the study of control systems is the notion of exponential stability (often referred to simply as stability).  
A related stronger condition known as superstability has been the focus of much research (e.g. \cite{Bal04}, \cite{Bal81}, \cite{NR92}, \cite{RW95}, \cite{Lu01}), in particular its connection with finite time extinction. The importance of exponentially stable semigroups stems, in part, from the fact that such systems admit transfer functions. Furthermore, if the semigroup is superstable (i.e $\omega_0=-\infty$) then $G(\lambda)$ can be defined on the whole complex plane.

The notion of superstability first appeared (in a very rough form) in the seminal work of Hille and Phillips \cite{HP57} who were concerned primarily with the mathematical aspects of it and related properties, particularly the relationship between the spectrum of the infinitesimal generator and the stability of the semigroup.
%In Hilbert spaces if the system $\Sigma (A,I,I)$ admits a transfer function \cite{Gearhart} then the semigroup generated by $A$ is exponentially stable, this doesn't hold in the general case of Banach spaces.
Later work refined and extended this notion, applying more complicated machinery (\cite{NR92}, \cite{RW95}).  More recently Balakrishnan \cite{Bal04} (and others \cite{udwadia2005}, \cite{udwadia2012}, \cite{liu2012}, \cite{shang2013}) have become interested in the superstability phenomena arising in the control theory of (physical) systems.

Particularly interesting for control applications is the relationship between superstability and finite time extinction (as opposed to asymptotic stability). 
In \cite{Bal04}, Balakrishnan poses the following question: are there physical (i.e. differential operator) examples of superstability which are not of extinction-in-finite-time type?  In this paper, we provide a positive answer to this question with a constructive example (in fact many similar examples can easily be obtained).  We also collect, clarify and extend the existing results in the field with a new unified approach. In practice, our results provide a clear methodology to check if a superstable system also exhibits finite time extinction. Thus our results also provide guidelines for the design of controllers to guarantee finite-time extinction.

The key ingredient in our work is a new approach to the concepts of stability and superstability focusing on the ``entry times'' of the system into balls about the origin (in the Banach space).  This approach, which has a certain probabilistic flavor though is not in itself probabilistic, allows us to unify the existing results in the field with (largely) new proofs.  More importantly, we obtain analogues of certain well-known results about stability for superstable and finite-time-extinction systems.  In particular, an analogue of Pazy's condition \cite{Paz83} for stability is given for both superstability and finite time extinction.

\section{Main Results}

Our main result characterizing the types of stability is:
\begin{theorem*}
Let $\{T(t)\}$ be a $C_{0}$-semigroup of bounded linear operators on a Banach space $X$.  Define the relative entry time for each $r \in \mathbb{N}$ as
\[
u_{r} = \sup \{ t_{r+1}(x) - t_{r}(x) : \| x \| \leq 1 \} \quad\text{where}\quad t_{r}(x) = \inf \{ t \geq 0 : \| T(t^{\prime}) x \| \leq e^{-r} \text{ for all $t^{\prime} \geq t$} \}.
\]
Then
\begin{align*}
\text{(i)}\quad & \text{$\{ T(t) \}$ is stable} & \text{if and only if} &&& \limsup_{r\to\infty} u_{r} < \infty; \\
\text{(ii)}\quad & \text{$\{ T(t) \}$ is superstable} & \text{if and only if} &&& \lim u_{r} = 0; \\
\text{(iii)}\quad & \text{$\{ T(t) \}$ has finite time extinction} & \text{if and only if} &&& \sum_{r} u_{r} < \infty.
\end{align*}
\end{theorem*}
Our Pazy-type condition is:
\begin{theorem*}
Let $\{ T(t) \}$ be a $C_{0}$-semigroup.  Then, if for some $a > 0$,
\begin{align*}
\text{(i)} \quad &\int_{a}^{\infty} \| T(t) \|^{p} dt < \infty\quad \text{for some $0 < p < \infty$} \quad & \text{then} \quad & \text{$\{T(t)\}$ is stable}; \\
\text{(ii)} \quad &\int_{a}^{\infty} \big{|} \log \| T(t) \| \big{|}^{-p} dt < \infty \quad \text{for some $1 < p < \infty$} \quad & \text{then} \quad & \text{$\{T(t)\}$ is stable}; \\
\text{(iii)} \quad &\int_{a}^{\infty} \big{|}\log \| T(t) \| \big{|}^{-1} dt < \infty \quad & \text{then} \quad & \text{$\{T(t)\}$ is superstable}; \\
\text{(iv)} \quad &\lim_{p\downarrow 0} \int_{a}^{\infty} \big{|} \log \| T(t) \| \big{|}^{-p} dt < \infty \quad & \text{then} \quad & \text{$\{T(t)\}$ has finite time extinction}.
\end{align*}
\end{theorem*}

\section{Preliminaries}
A family $\{T(t)\}_{t\geq 0}$ of bounded linear operators on a Banach space $X$ is called a
 {\bf strongly continuous semigroup} (or {\bf $C_{0}$-semigroup})
 when $T(0)=I, \; T(t+s)=T(t)T(s)$ for all $t,s \geq 0$, and
 $\lim _{t \downarrow 0} T(t) = I$ in the strong operator topology ($T(t)x \to x$ as $t \downarrow 0$
for all $x \in X$).
As is well-known, this implies that the map $t \mapsto T(t)$ is (strongly) continuous.

For a strongly continuous semigroup $\{T(t)\}_{t\ge 0}$, define $\mathcal{D}$ to be the set of all $x\in X$ such
that $\lim_{t\downarrow 0}t^{-1}(T(t)x-x)$ exists. The {\bf
infinitesimal generator} of the semigroup $\{T(t)\}_{t\ge
0}$ is the operator $A$ on $X$, with the domain $D(A)=
\mathcal{D}$, given by
\[
Ax=\lim\limits_{t\downarrow 0}\frac{T(t)x-x}{t},\;x\in D(A).
\]
The name ``infinitesimal generator'' is justified by the fact that
\[
Ax=\frac{d(T(t)x)}{dt}\Big{|}_{t=0},\;x\in D(A).
\]
The pair $(A,\mathcal{D})$ and the semigroup $\{T_{t}\}_{t\geq 0}$ uniquely determine one another.

A strongly continuous semigroup $\{T(t)\}_{t\geq 0}$ is called {\bf exponentially stable} (or just {\bf stable}) when there exists constants $M > 0$ and $\rho > 0$ such that
\[
\| T(t) \| \leq M e^{-\rho t} \quad\text{for all $t \geq 0$}.
\]
Equivalently, define
the {\bf stability index} to be
\[
\sup \{ \nu > 0 : (\exists M > 0) \| T(t) \| \leq M e^{-\nu t} \quad\text{for all $t \geq 0$} \},
\]
and stability is then the requirement that the index be positive.

The {\bf growth characteristic} is
\[
\omega_{0} = \lim_{t\downarrow 0} \frac{\log \| T(t) \|}{t}
\]
which, as is well-known, is equal to $-\nu$ where $\nu$ is the stability index (when the semigroup is stable).

It is then natural to define {\bf superstability} to be the condition that the growth characteristic is $\omega_{0} = -\infty$.  Alternatively, superstability can be defined as the equivalent condition that the operators $T(t)$ be quasinilpotent (recall that an operator $T$ is quasinilpotent when $spec(T) = \{ 0 \}$).

A system is said to have {\bf finite time extinction} when there is some $t_{0} \geq 0$ such that $T(t)x = 0$ for all $t \geq t_{0}$ and all $x \in X$ with $\| x \| \leq 1$.

A semigroup $\{ T(t) \}$ is {\bf nilpotent} when there exists $t_{0}$ such that $T(t_{0}) = 0$.  The smallest possible choice $t_{0}$ such that $T(t^{\prime}) = 0$ for all $t^{\prime} > t_{0}$ is called the {\bf index of nilpotency} for the semigroup.

The reader should note that in what follows we often defer the proofs until after all results are stated, the purpose being to stress the similarities among the theorems characterizing these three concepts (perhaps the most useful aspect of our approach).

\section{Final Entry Times}

\begin{definition}
Let $\{T(t)\}$ be a $C_{0}$-semigroup on a Banach space $X$.  For each $x \in X$ and $r \in \mathbb{N}$ the \textbf{\emph{final entry time}} of $x$ into the $e^{-r}$-ball is
\[
t_{r}(x) := \inf \{ t \geq 0 : \| T(t^{\prime}) x \| \leq e^{-r} \text{ for all $t^{\prime} \geq t$} \},
\]
and the final entry time of the $1$-ball into the $e^{-r}$-ball (referred to from here on as just the final entry time of the $e^{-r}$-ball) is
\[
t_{r} := \sup \{ t_{r}(x) : \| x \| \leq 1 \},
\]
where we adopt the (usual) convention that the infimum of the empty set is $\infty$ (that is, $t_{r}(x) = \infty$ when no such time exists).\\
\\
The \textbf{\emph{relative entry time}} of $x$ into the $e^{-r}$-ball is
\[
u_{r}(x) := \begin{cases} t_{r+1}(x) - t_{r}(x) &\quad\text{when $t_{r+1}(x) < \infty$} \\ \infty &\quad\text{when $t_{r+1}(x) = \infty$} \end{cases}
\]
and the \textbf{\emph{relative entry time}} of the $1$-ball into the $e^{-r}$-ball (referred to from here on as just the final entry time of the $e^{-r}$-ball) is
\[
u_{r} := \sup \{ u_{r}(x) : \| x \| \leq 1 \}.
\]
\end{definition}

\begin{lemma}\label{L:2}
For any $r$ we have $u_{r} = t_{r+1} - t_{r}$ (meaning when either or both $t_{r}, t_{r+1} = \infty$ then $u_{r} = \infty$).
\end{lemma}
\begin{proof}
First note that $u_{r} = \infty$ if and only if $t_{r+1} = \infty$ and that $t_{r} = \infty$ implies $t_{r+1} = \infty$ so we need only handle the case that all three are finite.
By definition,
\[
u_{r} + t_{r} = \sup \{ u_{r}(x) : \|x\|\leq 1 \} + \sup \{ t_{r}(x) : \|x\|\leq 1\} \geq t_{r+1}.
\]
On the other hand, there exists a sequence $x_{n}$ such that $t_{r}(x_{n}) \to t_{r}$.  So for any $\epsilon > 0$ and sufficiently large $n$ we have $t_{r}(x_{n}) > t_{r} - \epsilon$ and so $t_{r} < t_{r}(x_{n}) + \epsilon$.  By definition $t_{r+1} \geq t_{r+1}(x)$ for all $x$, in particular for the $x_{n}$.
Then 
\[
t_{r+1} - t_{r} \geq t_{r+1}(x_{n}) - t_{r}(x_{n}) - \epsilon \geq u_{r} - \epsilon.
\]
As $\epsilon$ is arbitrary, the claim follows.
\end{proof}

\begin{lemma}\label{L:e-r}
For any $x$ and $r$ with $\| x \| > e^{-r}$ we have $\| T(t_{r}(x)) x \| = e^{-r}$.  Moreover, $\| T(t_{r}) \| = e^{-r}$.
\end{lemma}
\begin{proof}
This follows directly from the strong continuity of $T(t)$ (which is automatic for right continuity and follows from the uniform boundedness principle for left continuity).
\end{proof}

\begin{lemma}\label{L:1}
The sequence $u_{r}$ is nonincreasing in $r$: $u_{r+1} \leq u_{r}$.  Hence the limit $\lim_{r\to\infty} u_{r}$ always exists.
\end{lemma}
\begin{proof}
First note that if $u_{r}$ ever attains $\infty$ then in fact $t_{r+1} = \infty$ at that point and for all subsequent $t_{r}$ meaning that $u_{r}$ remains $\infty$ from then on.  So we need only handle the case when $u_{r} < \infty$ for all $r$ (and therefore assume that $t_{r} < \infty$ for all $r$). 

Now suppose that $u_{r} > u_{r-1}$ for some $r$.  Then there is some $x$ such that $u_{r}(x) > u_{r-1}$ (since $u_{r}$ is the supremum over all $x$).
Set 
\[
y = e T(t_{r}(x)) x \quad\text{ so that }\quad T(t)y = e T(t + t_{r}(x))x \quad\text{ for all $t \geq 0$}.
\]
By definition of $t_{r+1}(x)$ we have that
\[
\| T(t + t_{r}(x))x \| \leq e^{-r-1} \quad\text{ if and only if }\quad t + t_{r}(x) \geq t_{r+1}(x).
\]
Then 
\[
\| T(t) y \| = e \| T(t + t_{r}(x)) x \| \leq e^{-r} \quad\text{ if and only if }\quad t \geq t_{r+1}(x) - t_{r}(x) = u_{r}(x),
\]
which means that $t_{r}(y) = u_{r}(x)$.
Now $t_{r-1}(y) = 0$ since $\| T(t^{\prime} + t_{r}(x)) x \| \leq e^{-r}$ by definition of $t_{r}(x)$ and so
\[
\| T(t^{\prime}) y \| = e \| T(t^{\prime} + t_{r}(x)) x \| \leq e e^{-r} = e^{-(r-1)} \quad\quad\text{ for all $t^{\prime} \geq 0$}.
\]
Hence $u_{r-1}(y) = u_{r}(x)$.  But this means that
\[
u_{r}(x) > u_{r-1} \geq u_{r-1}(y) = u_{r}(x),
\]
which is a contradiction.  Therefore $u_{r} \leq u_{r-1}$ for all $r$ as claimed.
\end{proof}

\begin{lemma}\label{L:stoptime}
If the $T(t)$ are (not necessarily proper) contractions (i.e., $\| T(t) \| \leq 1$) then
\[
t_{r}(x) = \sup \{ t \geq 0 : \| T(t) x \| \geq e^{-r} \},
\]
which is to say the $t_{r}$ are ``stopping times''.
\end{lemma}
\begin{proof}
The fact that $T(t)$ are contractions forces that for any $q > 0$: 
\[
\| T(t+q)x \| = \| T(q) T(t) x \| \leq \| T(t)x \| \| T(q) \| \leq \| T(t) x \|,
\]
by the semigroup property.
\end{proof}

\section{The Entry Time Growth Characteristic $\omega_{0}^{(ET)}$}

\begin{definition}
The \textbf{\emph{entry time growth characteristic}} of a $C_{0}$-semigroup $\{ T(t) \}$ is defined by
\[
\omega_{0}^{(ET)} := -\big{(} \lim_{r\to\infty} u_{r} \big{)}^{-1},
\]
which always exists by Lemma \ref{L:1}.
\end{definition}

We now show that the entry time growth characteristic is equal to the usual growth characteristic $\omega_{0}$ defined previously for stable semigroups.  Note that if the semigroup is not stable then $\omega_{0}^{(ET)} = 0$.

\begin{theorem}\label{T:omega0}
For a stable $C_{0}$ semigroup $T(t)$,
\[
\omega_{0}^{(ET)} = \omega_{0} = \lim_{t\downarrow 0} \frac{\log \| T(t) \| }{t} = \inf_{t\geq 0}\frac{\log \| T(t) \|}{t} = \lim_{t\to\infty} \frac{\log \| T(t) \|}{t}.
\]
\end{theorem}
\begin{proof}
That $\omega_{0}$ (the usual growth characteristic), which is defined as the limit when $t$ goes down to zero above, is equal to the infimum is well-known.

Define $\omega(t) = \frac{1}{t}\log \| T(t) \|$.  Then,
\[
(t + s) \omega(t + s) = \log \| T(t + s) \| \leq \log \| T(t) \| \| T(s) \|
= \log \| T(t) \| + \log \| T(s) \|
= t \omega(t) + s \omega(s)
\]
is a subadditive sequence and therefore
\[
\lim_{t \to \infty} \frac{t \omega(t)}{t} = \lim_{t \to \infty} \omega(t)
\]
exists.  Call this limit $\omega$.  Set $w_{r} = \omega(t_{r})$.  By Lemma \ref{L:e-r} we know that $w_{r} = - \frac{r}{t_{r}}$.
Now $\lim_{r\to\infty} w_{r} = \omega$ and
\[
- \frac{r+1}{t_{r+1}} = - \frac{r}{t_{r}}\frac{t_{r}}{t_{r+1}} - \frac{1}{t_{r+1}},
\]
and so
\[
w_{r+1} = w_{r} \frac{t_{r}}{t_{r+1}} - \frac{1}{t_{r+1}},
\]
which means that
\[
w_{r+1} t_{r+1} - w_{r} t_{r} = -1.
\]
Consider the case when $\omega > -\infty$.  Taking limits in the above,
$\lim_{r\to\infty} \omega u_{r} = -1$
and so $\omega_{0}^{(ET)} = \omega$.

Now consider the case when $\omega = -\infty$.  Suppose that $\omega_{0}^{(ET)} > -\infty$.  Then $\lim_{r\to\infty} u_{r} = c > 0$ and so
\[
-1 = t_{r+1}w_{r+1} - t_{r} w_{r} = t_{r} (w_{r+1} - w_{r}) + u_{r} w_{r+1} \leq u_{r}w_{r+1},
\]
since $w_{r+1} - w_{r} \leq 0$.  But then
\[
-1 \leq \lim_{r\to\infty} u_{r} w_{r+1} = c (-\infty) = -\infty,
\]
contradicting that $\omega_{0}^{(ET)} > -\infty$.
Therefore in either case $\omega_{0}^{(ET)} = \omega$.

Now we show that $\omega = \omega_{0}$.  Since we know $\omega_{0} = \inf_{t \geq 0} \omega(t)$ and $\omega = \lim_{t \to \infty} \omega(t)$ we already have that $\omega_{0} \leq \omega$.  
For any positive integer $n$ and any $s \geq 0$ we have
\[
\omega(ns) = \frac{1}{ns} \| T(s)^{n} \| \leq \frac{1}{ns} \log \| T(s) \|^{n} = \omega(s),
\]
and therefore
\[
\lim_{t \to \infty} \omega(t) = \lim_{n \to \infty} \omega(ns) \leq \omega(s),
\]
which means that
\[
\omega = \lim_{t \to \infty} \omega(t) \leq \inf_{t \geq 0} \omega(t) = \omega_{0}.
\]
Therefore $\omega = \omega_{0}$ and the proof is complete.
\end{proof}

\section{Equivalence of Stability Notions}

We now present the main theorems characterizing the various notions of stability.  One of our aims in this paper is to collect and clarify the various characterizations of these notions.  To this end
we include several known characterizations and provide new proofs using our techniques.
Specifically, the equivalences in this section, excepting the conditions which involve the relative entry times $u_{r}$, are known (see e.g. \cite{Bal04}).

\begin{theorem}\label{T:stable}
For a $C_{0}$-semigroup $\{ T(t) \}$ and $\nu > 0$, the following are equivalent:
\begin{align*}
\text{(i)} \quad&\text{$\{T(t)\}$ is \textbf{\emph{stable}} with \textbf{\emph{stability index}} $\nu$}: \quad(\forall \rho < \nu) (\exists M > 0) (\forall t \geq 0) \| T(t) \| \leq Me^{-\rho t}; \\
\text{(ii)}  \quad&\lim_{r\to\infty} u_{r} = \nu^{-1} < \infty; \\
\text{(iii)}  \quad&spec (T(t)) \subseteq \{ \lambda \in \mathbb{C} : |\lambda| \leq e^{-t\nu} \}; \\
\text{(iv)}  \quad&\omega_{0} = -\nu.
\end{align*}
\end{theorem}

\begin{theorem}\label{T:stableA}
If $\{ T(t) \}$ is a stable $C_{0}$-semigroup with stability index $\nu$ and $A$ is the generator of $T(t)$ with domain $D$ then $spec (A) \subseteq \{ \lambda \in \mathbb{C} : Re\ \lambda < -\nu \}$.  The converse does not hold.
\end{theorem}

\begin{remark}
The spectrum of the generator, $spec (A)$, depends very delicately on the domain of definition (see the counterexamples below).  In particular, the operator $A$ treated as having full domain may well have spectrum much larger than $spec (A)$.
\end{remark}

\begin{theorem}\label{T:superstable}
For a $C_{0}$-semigroup $\{ T(t) \}$, the following are equivalent:
\begin{align*}
\text{(i)}  \quad&\text{$\{T(t)\}$ is \textbf{\emph{superstable}}}: \quad\text{$\{T(t)\}$ is stable with stability index $\infty$}; \\
\text{(ii)}  \quad&\lim_{r\to\infty} u_{r} = 0; \\
\text{(iii)}  \quad&(\forall \nu > 0) (\exists M_{\nu} > 0) (\forall t \geq 0) \| T(t) \| \leq M_{\nu}e^{-\nu t}; \\
\text{(iv)}  \quad&\text{$T(t)$ are all quasinilpotent: $spec (T(t)) = \{ 0 \}$}; \\
\text{(v)}  \quad&\omega_{0} = -\infty. \\
\end{align*}
\end{theorem}

\begin{remark}
The constants $M_{\nu}$ in condition (iii) must tend to infinity as $\nu \to \infty$ since otherwise the semigroup will be identically $0$.
\end{remark}

\begin{theorem}\label{T:superstableA}
If $\{ T(t) \}$ is a superstable $C_{0}$-semigroup with generator $A$ and domain $D$ then $spec (A) = \emptyset$.  The converse does not hold.
\end{theorem}

\begin{theorem}\label{T:finitetime}
For a $C_{0}$-semigroup $\{ T(t) \}$ on a Banach space $X$ with generator $A$ having domain $D$, and $0\leq k < \infty$ the following are equivalent:
\begin{align*}
\text{(i)}  \quad&\text{ $\{T(t)\}$ has \textbf{\emph{finite time extinction at time $k$}}}: \\
&\quad\quad(\forall x \in X) (\exists t_{\infty}(x) \geq 0) (\forall t \geq t_{\infty}(x)) T(t)x = 0 \text{ and } \sup 
\{t_{\infty}(x) : \|x\| \leq 1\} = k; \\
\text{(ii)}  \quad&\sum_{r} u_{r} = k < \infty; \\
\text{(iii)}  \quad&(\forall \nu > 0) (\exists M_{\nu} > 0) (\forall t \geq 0) \| T(t) \| \leq M_{\nu} e^{-\nu t} \text{ and } \sup_{\nu > 0} \frac{\log M_{\nu}}{\nu} = k; \\
\text{(iv)}  \quad&(\exists M > 0) (\forall \nu \geq 0) (\forall t \geq 0) \| T(t) \| \leq M e^{- \nu (t - k)}; \\
\text{(v)}  \quad&\text{$T(t)$ is nilpotent with nilpotency index $k$: $T(q) = 0$ for $q > k$ and $T(q) \ne 0$ for $q < k$}; \\
\text{(vi)}  \quad&\text{the resolvent function $R(\lambda,A)$ is entire } \\
&\quad\quad\text{ and } \big{|}R(\lambda,A)\big{|} \leq C(1 + |\lambda|)^{-N}e^{k |Re \lambda|} \text{ for some constants $C,N$}. \\
\end{align*}
\end{theorem}

\begin{remark}
In condition (i), the definition of finite time extinction, the  $t_{\infty}$ can be chosen uniformly on bounded sets, in particular on balls around the origin of finite radius, however, $t_{\infty}$ cannot be chosen uniformly over $x$ unless the underlying space of the Banach space is compact (i.e., $L^{2}[0,1]$ not $L^{2}[0,\infty)$).
\end{remark}

\begin{theorem}
Extinction in finite time implies superstability and superstability implies stability.  The converses of both statements are false.
\end{theorem}

\begin{theorem}\label{T:C}
Let $\{ T_{t} \}$ be a $C_{0}$-semigroup.  Then, if for some $a > 0$,
\begin{align*}
\text{(i)} \quad &\int_{a}^{\infty} \| T(t) \|^{p} dt < \infty\quad \text{for some $0 < p < \infty$} \quad & \text{then} \quad & \text{$\{T(t)\}$ is stable}; \\
\text{(ii)} \quad &\int_{a}^{\infty} \big{|} \log \| T(t) \| \big{|}^{-p} dt < \infty \quad \text{for some $1 < p < \infty$} \quad & \text{then} \quad & \text{$\{T(t)\}$ is stable}; \\
\text{(iii)} \quad &\int_{a}^{\infty} \big{|}\log \| T(t) \| \big{|}^{-1} dt < \infty \quad & \text{then} \quad & \text{$\{T(t)\}$ is superstable}; \\
\text{(iv)} \quad &\lim_{p\downarrow 0} \int_{a}^{\infty} \big{|} \log \| T(t) \| \big{|}^{-p} dt < \infty \quad & \text{then} \quad & \text{$\{T(t)\}$ has finite time extinction}.
\end{align*}
\end{theorem}

\begin{remark}
Stability can only occur when $\| T(t) \|$ is eventually bounded by $1$ (that is, $t_{0} < \infty$).  Taking $a = t_{0}$ will cause the integrals to converge whenever there is some value of $a$ that causes convergence.
\end{remark}

\section{Proofs of Equivalences}

As usual, let $\{T(t)\}$ be a $C_{0}$-semigroup with generator $A$ having dense domain $D$ on the Banach space $X$.

\begin{proof} (of Theorem \ref{T:stable}).
Assume condition (ii) holds:
\[
\lim_{r\to\infty} u_{r} = \nu^{-1} < \infty.
\]
  Fix $\epsilon > 0$.  Then there exists $r^{*}$ such that
  \[
  u_{r} < \nu^{-1} + \epsilon \quad\text{ for all }\quad r \geq r^{*}.
  \]
    Hence for $r \geq r^{*}$, 
    \[
    t_{r} - t_{r^{*}} \leq (r - r^{*})(\nu^{-1} + \epsilon).
    \]
For any $t \geq r^{*} (\nu^{-1} + \epsilon)$ pick $r \geq r^{*}$ such that
\[
r(\nu^{-1} + \epsilon) + t_{r^{*}} \leq t < (r + 1)(\nu^{-1} + \epsilon) + t_{r^{*}}.
\]
Since $t \geq t_{r}$ (as $r(\nu^{-1} + \epsilon > 0$) we have that $\| T(t) \| \leq e^{-r}$ and 
\[
t < (r + 1)(\nu^{-1} + \epsilon) + t_{r^{*}} \quad\text{implies}\quad r > \frac{t - t_{r^{*}}}{\nu^{-1} + \epsilon} - 1.
\]
Then
\[
\| T(t) \| \leq e^{-r} < e^{1 - \frac{t - t_{r^{*}}}{\nu^{-1} + \epsilon}} = e^{1 + \frac{t_{r^{*}}}{\nu^{-1} + \epsilon}} e^{-t \frac{1}{\nu^{-1}+\epsilon}}.
\]
So $\{ T(t) \}$ has stability index less than $\frac{1}{\nu^{-1} + \epsilon}$.  Since $\epsilon$ was arbitrary, condition (i) holds.

Conversely, assume (i) holds.  Then for any $\rho < \nu$ there is $M$ such that $\|T(t)\| \leq Me^{-\rho t}$.  For $t \geq \frac{r - \log M}{\rho}$ we then have that $\| T(t) \| \leq e^{-r}$.  Hence
\[
t_{r} \leq \frac{r - \log M}{\rho}.
\]
Suppose $\lim u_{r} > \frac{1}{\rho} + 2\delta$ for some $\delta > 0$.  Then for sufficiently large $r^{\prime}$ we have $u_{r} > \frac{1}{\rho} + \delta$ for $r \geq r^{\prime}$ so
\[
t_{r + 1} = t_{r^{\prime}} + u_{r^{\prime}} + \cdots + u_{r} > t_{r^{\prime}} + (r - r^{\prime})\frac{1}{\rho} + (r - r^{\prime})\delta,
\]
and therefore
\[
\frac{r + 1 - \log M}{\rho} - (r - r^{\prime})\frac{1}{\rho} > (r - r^{\prime})\delta,
\]
so
\[
\frac{r^{\prime} + 1 - \log M}{\rho} + r^{\prime}\delta > r \delta,
\]
but the left-hand side is constant and the right hand tends to $\infty$ as $r \to \infty$.  This
means that $\lim u_{r} \leq \frac{1}{\rho}$.  Since $\rho < \nu$ is arbitrary we have (ii).

Now assume (iv) holds.  Then by Gelfand's spectral radius formula,
\[
\sup \big{|} spec (T(t)) \big{|} = \lim_{t} \| T(t) \|^{\frac{1}{t}} = e^{\omega_{0} t} = e^{-t \nu}. 
\]
Hence (iii) holds.  Likewise, if (iii) holds then by Gelfand's formula, $\omega_{0} = -\nu$ so (iv) holds.

The equivalence of (i) and (ii) with (iv) is a direct consequence of Theorem \ref{T:omega0}.  This completes the proof that (i) through (iv) are equivalent.
\end{proof}

\begin{proof} (of Theorem \ref{T:stableA}).
This is well-known.
\end{proof}

\begin{proof} (of Theorem \ref{T:superstable}).
The equivalences follow from identical arguments to those for the case of stable semigroups (simply replace $\nu^{-1}$ by $0$).  
\end{proof}

\begin{proof} (of Theorem \ref{T:superstableA}).
This follows from Theorem \ref{T:stableA}: if $z \in spec(A)$ then the stability index $\nu$ is at least $- Re(z)$ but a superstable semigroup has stability index $\infty$.
\end{proof}

\begin{proof} (of Theorem \ref{T:finitetime}).
Assume that (ii) holds.  Then for any $x$ with $\| x \| \leq 1$ we have that 
\[
t_{r+1}(x) - t_{0}(x) = \sum_{j=0}^{r} u_{r}(x) \leq \sum_{j=0}^{r} u_{r} \to k,
\]
hence (i) holds.

Conversely, assume (i) holds and suppose (ii) fails.  If $\sum_{r} u_{r} = \ell < k$ then
\[
t_{\infty}(x) \leq \ell < k \quad\text{for all $x$},
\]
contradicting (i).  So it must be that $\sum_{r} u_{r} = \ell > k$.  By Lemma \ref{L:2} we have
\[
\sum_{r}u_{r} = \sum_{r} t_{r+1} - t_{r} = t_{\infty} - t_{0},
\]
and there is then a sequence $x_{n}$ such that $t_{\infty}(x_{n}) \to \ell > k$ contradicting (i).

Assume (ii) holds.  For $\nu > 0$ pick $r^{*}$ such that
\[
\sup_{r\geq r^{*}} u_{r} < \nu^{-1} \quad\text{and set}\quad M_{\nu} = e^{1 + \nu t_{r^{*}}}.
\]
Then, as in the proof of stability, we have
\[
r\nu^{-1} \leq t - t_{r^{*}} < (r+1)\nu^{-1}\quad\text{implies}\quad
\|T(t)\| \leq M_{\nu} e^{-t \nu}.
\]
Now
\[
\frac{\log M_{\nu}}{\nu} = \frac{1}{\nu} + t_{r^{*}(\nu)}
\quad
\text{and}
\quad
r^{*}(\nu) \to \infty \quad\text{as $\nu \to \infty$},
\]
hence $\frac{\log M_{\nu}}{\nu} \to \sum_{r} u_{r} = k$.  So (iii) holds.

Assume (iii) holds.  Then
\[
M_{\nu} \leq e^{k\nu}\quad\text{for all $\nu$},
\]
hence (iv) holds with $M = 1$.

Assume (iv) holds.  Then for $t > k$ we have
\[
\| T(t) \| \leq Me^{-\nu(t-k)} \quad\text{for all $\nu$ and $t - k > 0$},
\]
so $\lim_{\nu} e^{-\nu(t-k)} = 0$ hence $\| T(t) \| = 0$.  Thus (i) holds.

The equivalence of (i) and (v) is trivial.

To see that (i) and (vi) are equivalent, note that (i) and (vi) both imply $spec(A) = \emptyset$.
Then
\[
R(\lambda, A) = \int_{0}^{\infty} e^{-\lambda t}T(t) dt \quad\text{for all $\lambda$}
\]
(in general for $Re\ \lambda > \omega_{0}$, see, e.g,. \cite{Bal81}).

Hence $R(i\lambda,A)$ is the Fourier Transform of $T(t)$.  By the Paley-Wiener Theorem, $R(i\lambda,A)$ is the Fourier Transform of a compactly supported function (i.e. $T(t) = 0$ for all $t \geq k$) if and only if (vi) holds.  This argument first appeared in \cite{GK70}.
\end{proof}

\section{Proof of the Pazy-type Criteria}

\begin{lemma}\label{L:Ftrick}
Suppose that $\| T(t) \| \leq 1$ for all $t \geq 0$ (that is, $t_{0} = 0$).
Let $F: \mathbb{R} \cup \{ \pm \infty \} \to [0,\infty]$ be a decreasing function such that $F(\infty) = 0$.  Then
\[
\sum_{r=0}^{\infty} u_{r} F(r+1) \leq \int_{0}^{\infty} F(-\log \|T(t)\|) dt  \leq \sum_{r=0}^{\infty} u_{r} F(r).
\]
\end{lemma}
\begin{proof}
Observe that
\[
\int_{0}^{\infty} F(-\log \|T(t)\|) dt = \int_{0}^{t_{0}} F(-\log \|T(t)\|) dt + \sum_{r=0}^{\infty} \int_{t_{r}}^{t_{r+1}} F(-\log \|T(t)\|) dt + \int_{t_{\infty}}^{\infty} F(-\log \|T(t)\|) dt.
\]

Now $t_{0} = 0$ so the first term on the right is $0$.  For $t \geq t_{\infty}$ we have that $\|T(t)\| = 0$ so $F(-\log \|T(t)\|) = F(\infty) = 0$ meaning that the third term on the right is zero.

For the middle terms, note that for $t < t^{\prime}$,
\[
\| T(t^{\prime}) \| = \| T(t) T(t^{\prime} - t) \| \leq \| T(t) \| \| T(t^{\prime} - t) \| \leq \| T(t) \|,
\]
since $\| T(t^{\prime} - t) \| \leq 1$ by assumption.  So
\[
\| T(t_{r+1}) \| \leq \| T(t) \| \leq \| T(t_{r}) \| \quad\text{for}\quad t_{r} \leq t \leq t_{r+1},
\]
and therefore, by Lemma \ref{L:e-r},
\[
-r-1 = \log \| T(t_{r+1}) \| \leq \log \| T(t) \| \leq \log \| T(t_{r}) \| = -r \quad\text{for}\quad t_{r} \leq t \leq t_{r+1},
\]
and so since $F$ is decreasing
\[
F(r+1) \leq F(-\log \|T(t)\|) \leq F(r) \quad\text{for}\quad t_{r} \leq t \leq t_{r+1}.
\]
So for each term in the sum,
\[
u_{r}F(r+1) = \int_{t_{r}}^{t_{r+1}} F(r+1) dt \leq \int_{t_{r}}^{t_{r+1}} F(-\log \|T(t)\|) dt \leq u_{r} F(r).
\]
\end{proof}

\begin{proof} (of Theorem \ref{T:C}).
First note that if $\limsup \| T(t) \| > 1$ then $u_{r} = \infty$ for all $r$ so there can be no stability.  In this case, none of the three conditions involving integrals can hold.  So it is enough to consider the case when $\| T(t) \|$ is eventually bounded by $1$.  Since
\[
\int_{0}^{t_{0}} H(\| T(t) \|) dt < \infty
\]
for any bounded function $H$ (recall that $t_{0} < \infty$ since we have eliminated the other case), we may assume that $\| T(t) \| \leq 1$ for all $t$: the integral conditions are unaffected by finite translations in time as is the stability of the semigroup.

Recall that condition (i) is the Datko-Pazy Theorem (\cite{Dat70}, \cite{Paz72}, \cite{Paz83}).  Consider the function $F(x) = e^{-px}$ for some fixed $0 < p < \infty$.  Then $F$ is decreasing and $F(\infty) = 0$.

By Lemma \ref{L:Ftrick},
\[
\sum_{r=0}^{\infty} u_{r} e^{-p(r+1)} \leq \int_{0}^{\infty} F(-\log \|T(t)\|) dt = \int_{0}^{\infty} \|T(t)\|^{p} dt < \infty.
\]
For any given $r^{*}$ observe that since the $u_{r}$ are nonincreasing,
\[
\sum_{r=0}^{r^{*}-1} u_{r} e^{-p(r+1)} \geq  u_{r^{*}} \sum_{r=1}^{r^{*}} e^{-pr},
\]
and since $\sum_{r} e^{-p(r+1)} = C < \infty$,
\[
\lim_{r^{*}} \sum_{r=0}^{r^{*}-1} u_{r} e^{-p(r+1)} \geq \lim_{r^{*}} u_{r^{*}} C,
\]
and therefore
\[
\lim_{r\to\infty} u_{r} \leq C^{-1} \sum_{r=0}^{\infty} u_{r} e^{-p(r+1)} < \infty,
\]
hence the semigroup is stable.

Condition (ii) is a weakening of the Pazy condition: set $F(x) = x^{-p}$ for the appropriate $1 < p < \infty$.  By Lemma \ref{L:Ftrick},
\[
\sum_{r=0}^{\infty} u_{r} (r+1)^{-p} \leq \int_{0}^{\infty} (- \log \| T(t) \|)^{-p} dt < \infty.
\]
Then, as above, since $\sum_{r} (r+1)^{-p} = C < \infty$,
\[
\lim_{r\to\infty} u_{r} \leq C^{-1} \sum_{r=0}^{\infty} u_{r} (r+1)^{-p} < \infty,
\]
so the semigroup is stable.

Now condition (iii): set $F(x) = \frac{1}{x}$.  By Lemma \ref{L:Ftrick},
\[
\sum_{r=0}^{\infty} u_{r} \frac{1}{r+1} \leq \int_{0}^{\infty} \frac{dt}{-\log \|T(t)\|} < \infty.
\]
Proceeding as above,
\[
\sum_{r=0}^{r^{*}-1} u_{r} \frac{1}{r+1} \geq u_{r^{*}} \sum_{r=1}^{r^{*}} \frac{1}{r},
\]
and since $\sum_{r=1}^{\infty} \frac{1}{r} = \infty$ this means that $\lim_{r\to\infty} u_{r} = 0$ (and in fact converges to $0$ faster than the inverse of the harmonic sum $(1 + \cdots + \frac{1}{r})^{-1}$).  The semigroup is therefore superstable.

Finally condition (iv): set $F_{p}(x) = x^{-p}$ for $0 < p$.  Then, by Lemma \ref{L:Ftrick},
\[
\sup_{p>0} \sum_{r=0}^{\infty} u_{r} (r+1)^{-p} \leq \sup_{p>0} \int_{0}^{\infty} (- \log \| T(t) \|)^{p} dt < \infty,
\]
here we use that $F_{p} \leq F_{p^{\prime}}$ for $p \geq p^{\prime}$ so the hypothesis for $p\downarrow 0$ in fact implies boundedness for all $p > 0$.

Suppose that $\sum_{r} u_{r} = \infty$.  Then for any $K$ there exists $r_{K}$ such that $\sum_{r=0}^{r_{K}-1} u_{r} \geq K$.  Then
\[
\sum_{r=0}^{\infty} u_{r} (r+1)^{-p} \geq \sum_{r=0}^{r_{K}-1} u_{r} r_{K}^{-p} \geq K r_{K}^{-p},
\]
for any $p > 0$ and so
\[
\sup_{p>0} \sum_{r=0}^{\infty} u_{r} (r+1)^{-p} \geq K.
\]
But $K$ is arbitrary so
\[
\infty > \sup_{p>0} \sum_{r=0}^{\infty} u_{r} (r+1)^{-p} = \infty
\]
is a contradiction.  The semigroup therefore has finite time extinction.  In fact the semigroup goes extinct at time
\[
k = \sum_{r=0}^{\infty} u_{r} = \sup_{p>0} \int_{0}^{\infty} (- \log \|T(t)\|)^{-p} dt
\]
(details here are left to the reader).
\end{proof}

\section{Counterexamples}
We construct examples of semigroups demonstrating that finite time extinction is strictly stronger than superstability and that superstability is strictly stronger than stability.  In particular, we answer a question of Balakrishnan \cite{Bal04} on the existence of superstable semigroups not vanishing in finite time with generator being a differential operator (what he terms a ``physical system'').  We also remark on a (previously known) example showing that the spectrum of the generator does not fully determine superstability.

\subsection{Superstable Without Finite Time Extinction}

Consider the Gaussian (probability) measure $\mu $ on $\mathbb{R}^{+} = [0,\infty)$ given by $d\mu(x) = \sqrt{\frac{2}{\pi}} \exp(-\frac{x^{2}}{2}) dx$.  Let $X = L^{2}(\mathbb{R}^{+}, \mu)$.  Define the semigroup
\[
T(t)f(s) = f(s-t) \text{  for $s \geq t$ } \quad \text{ and } \quad T(t)f(s) = 0 \text{  otherwise }
\]
on $X$.  The reader may verify that this in fact a semigroup with generator $A = - \frac{d}{ds}$ and domain the appropriate Sobolev space.

Now for $f \in L^{2}(\mathbb{R}^{+}, \mu)$ with $\|f\|=1$ we have
\begin{align*}
\| T(t) f \|^{2} &= \int_{0}^{\infty} |f(s-t)|^{2} d \mu(s) = \int_{0}^{\infty} |f(v)|^{2} \sqrt{\frac{2}{\pi}} e^{-\frac{(t+v)^{2}}{2}} dv \leq e^{-\frac{t^{2}}{2}} \int_{0}^{\infty} |f(v)|^{2} d \mu(v) = e^{-\frac{t^{2}}{2}},
\end{align*}
since $e^{-(t+v)^{2}} \leq e^{- t^{2}}$ for $t,v \geq 0$.  So $\|T(t)\| \leq \exp(-\frac{t^{2}}{2}) \to 0$ meaning that the semigroup is superstable.  However, $T(t) \ne 0$ for any $t$.

Taking $f$ to be a norm one (with respect to $\mu$) function concentrated near $0$ we see that $\| T(t) \| = \exp(-\frac{t^{2}}{4})$ and so $t_{r} = 2 \sqrt{r}$ and $u_{r} = 2 ( \sqrt{r+1} - \sqrt{r} ) \to 0$ but $\sum u_{r} = \infty$.

Hence superstability can occur without finite time extinction (even when the generator is merely a derivative).  The space $(\mathbb{R}^{+},\mu)$ is a variant of the classical Gaussian measure space which arises naturally in the context of stochastic systems and quantum systems, among other areas.  Our example can easily be extended to any system with a Gaussian measure (details are left to the interested reader).

\subsection{Finite Time Extinction}

Define the semigroup
\[
T(t)f(s) = f(s+t) \text{  for $s + t \leq 1$ } \quad \text{ and } \quad T(t)f(s) = 0 \text{  otherwise }
\]
on $X = \{ f \in L^{2}[0,1] : f(0) = f(1) = f^{\prime}(0) = f^{\prime}(1) = 0\}$.  The reader may verify that this is a semigroup with generator $A = \frac{d}{ds}$ and domain Sobolev space.  It is clear that $T(t) = 0$ for all $t \geq 1$ so this semigroup has finite time extinction.  In fact, $t_{r} = 1$ for all $r > 0$ so $u_{r} = 0$ for $r > 0$ and $\sum u_{r} = 1 < \infty$.

\subsection{Stable but Not Superstable}

For completeness, we mention the fairly trivial example $T(t)f(s) = e^{-\nu t}f(s)$ (so the generator is $A = \nu I$ and the domain is $\mathcal{D} = L^{2}$) is clearly stable with index $\nu$.  Here $t_{r}$ is defined by $e^{-\nu t_{r}} = e^{-r}$ so $t_{r} = r \nu^{-1}$ meaning $u_{r} = -\nu^{-1}$.

\subsection{Empty Spectrum (for the Generator) but not Superstable}

For completeness, we mention an example due to Hille and Phillips \cite{HP57} (chapter 23, section 16).

We first present a superstable semigroup which will be used to develop the actual example of interest.
Define
$T(t)f(s) := \frac{1}{\Gamma(t)}\int_{0}^{s}(s-u)^{t-1}f(u)du$.  That this is a semigroup follows from Euler integral identities.
The generator $A$ is the derivative of the convolution with $\log$ minus a constant:
$Af(s) = \frac{d}{ds}\int_{0}^{s}\log(s-u) f(u)du - \gamma f(s)$ ($\gamma$ is Euler's constant).
Then $spec(A) = \emptyset$ and $\| T(t) \| \approx \frac{1}{t\Gamma(t)}$.  So $\omega_{0} = -\infty$ and $T(t) \to 0$ but $T(t) \ne 0$ for any $t$.  This semigroup is in fact superstable but does not have finite time extinction.

With this construction in hand, we construct the desired example:
for $\xi \in \mathbb{C}$ with $Re\ \xi > 0$, define $J^{\xi}f(s) = \frac{1}{\Gamma(\xi)}\int_{0}^{s}(s-u)^{\xi - 1}f(u) du$.  When $\xi$ is taken to be a positive real this yields the semigroup above.  There is an analytic extension of $J^{\xi}$ to $\xi$ purely imaginary.  Let $T(t) = J^{it}$.  The generator of this semigroup is $iA$ where $A$ is the generator from above.  Then $spec(iA) = \emptyset$ but $0 \notin spec(T(t))$, e.g. the operators are not quasinilpotent hence not superstable.  

The reader is referred to \cite{HP57} for details on these semigroups.

\bibliography{Superstability}

}\end{document}